\documentclass[psamsfonts,reqno,12pt]{amsart}
\usepackage{graphicx}
\usepackage{graphics}
\usepackage{latexsym}
\usepackage{amsmath}
\usepackage{amssymb}
\usepackage[numbers,sort&compress]{natbib}

\newtheorem{theorem}{Theorem}[section]
\newtheorem{lemma}[theorem]{Lemma}
\newtheorem{proposition}[theorem]{Proposition}
\newtheorem{corollary}[theorem]{Corollary}

\theoremstyle{definition}
\newtheorem{definition}[theorem]{Definition}
\newtheorem{example}[theorem]{Example}

\theoremstyle{remark}

\newtheorem{problem}[theorem]{Problem}
\newtheorem{question}[theorem]{Question}
\numberwithin{equation}{section}

\def\sideremark#1{\ifvmode\leavevmode\fi\vadjust{\vbox
to0pt{\vss \hbox to 0pt{\hskip\hsize\hskip1em
\vbox{\hsize2cm\tiny\raggedright\pretolerance10000
\noindent#1\hfill}\hss}\vbox to8pt{\vfil}\vss}}}

\begin{document}

\title{GENERALIZED-LUSH SPACES AND \\ THE MAZUR-ULAM PROPERTY}

\author[Dongni Tan]{Dongni Tan}
\author[Xujian Huang]{Xujian Huang}
\author[Rui Liu]{Rui Liu}
\address[Dongni Tan]{Department of Mathematics, Tianjin University of Technology, Tianjin 300384, P.R. China}
\email{tandongni0608@gmail.com}
\address[Xujian Huang]{Department of Mathematics, Tianjin University of Technology, Tianjin 300384, P.R. China}
\email{Huangxujian86@gmail.com}

\address[Rui Liu]{Department of Mathematics and LPMC, Nankai
University, Tianjin 300071, P.R. China}
\email{ruiliu@nankai.edu.cn}



\thanks{D. Tan was supported by the the National Natural Science Foundation
of China (NSFC) 11201338, X. Huang was supported by the NSFC
11201339, and R. Liu was supported by the NSFC 11001134 and is the
corresponding author.}

\maketitle
\begin{abstract} We introduce a new class of Banach spaces, called generalized-lush spaces
(GL-spaces for short), which contains almost-CL-spaces, separable
lush spaces (specially, separable $C$-rich subspaces of $C(K)$), and
even the two-dimensional space with hexagonal norm. We
obtain that the space $C(K,E)$ of the vector-valued continuous functions is a
GL-space whenever $E$ is, and show that the GL-space is stable under $c_0$-, $l_1$- and $l_\infty$-sums.
As an application, we prove that the Mazur-Ulam property holds for a larger class of
Banach spaces, called local-GL-spaces, including all lush spaces and
GL-spaces. Furthermore, we generalize the stability properties of GL-spaces to local-GL-spaces. From this, we can obtain many examples of Banach spaces having the Mazur-Ulam property.

\end{abstract}
\section{Introduction}
The classical Mazur-Ulam theorem states that every surjective
isometry between normed spaces is a linear mapping up to
translation. In 1972, Mankiewiz \cite{Ma} extended this by showing
that every surjective isometry between the open connected subsets of
normed spaces can be extended to a surjective affine isometry on the
whole space. This result implies that the metric structure on the
unit ball of a real normed space constrains the linear structure of
the whole space. It is of interest to us whether this result can be
extended to unit spheres. In 1987, Tingley \cite{Ti} first studied
isometries on the unit sphere and raised the isometric extension
problem:
\begin{problem}\label{pr:1}
Let $E$ and $F$ be normed spaces with the unit spheres $S_E$ and
$S_F$, respectively. If $T : S_E\to S_F$ is a surjective isometry,
then does there exist a linear isometry $\widetilde{T} : E\to F$
such that $\widetilde{T}|_{S_E}=T$?
\end{problem}
\noindent There is a number of publications on this topic and many positive answers on special spaces, for example,
$l^{p}(\Gamma),\,L^p(\mu)$($0< p\leq \infty$), $C(K)$, even the
James spaces and the (modified) Tsirelson spaces (see
\cite{D8,DG07,liu07a,LZ,Ta1,Ta3,Ta5,Ta4} and the references
therein).

Recently in \cite{CD}, Cheng and Dong considered the extension
question of isometries between unit spheres of Banach space and
introduced the Mazur-Ulam property:
\begin{definition}
A Banach space $E$ is said to have the \emph{Mazur-Ulam property}
(briefly MUP) provided that for every Banach space $F$, every
surjective isometry $T$ between the two unit spheres of $E$ and $F$
is the restriction of a linear isometry between the two spaces.
\end{definition}
\noindent Cheng and Dong attacked the problem for the class of
CL-spaces admitting a smooth point and polyhedral spaces. Unfortunately their interesting attempt failed by a mistake at the
very end of the proof (also see the introduction in \cite{KM,TL}).
In \cite{KM}, Kadets and Mart\'{i}n proved that finite-dimensional
polyhedral Banach spaces have the MUP. Notice that the problem is
still open even in two dimension case. In \cite{TL}, the authors Tan
and Liu proved that every almost-CL-spaces admitting a smooth point
(specially, separable almost-CL-spaces) has the MUP.

Recall that R. Fullerton \cite{Fu} first introduced the notion of
CL-spaces. It was extended by Lima \cite{Lima1,Lima2} who introduced
the almost-CL-spaces and gave the examples of real CL-spaces which are
$L_1(\mu)$ and its isometric preduals, in particular $C(K)$,
where $C(K)$ is a compact Hausdorff space. The infinite-dimensional
complex $L_1(\mu)$ spaces are proved by Mart\'{i}n and Pay\'{a }\cite{M}
to be only almost-CL-spaces. Lush spaces were introduced
recently in \cite{BK} and have been extensively studied recently in
\cite{BK2,KM, BK3}. Such spaces are of importance to supply an
example of a Banach space $E$ with the numerical index
$n(E)<n(E^*)$. It thus gives a negative answer to a question which
has been latent since the beginning of the theory of numerical
indices in the seventies. Now, a natural and interesting question
is: ``\emph{Does every almost-CL-space, even every lush space, has
the MUP?}"


In this paper, we introduce a natural concept of generalized-lush
spaces (GL-spaces for short), which contains almost-CL-spaces,
separable lush spaces (specially, separable $C$-rich subspaces of
$C(K)$), and even the two-dimensional space with hexagonal norm. We
obtain that the space $C(K,E)$ of the vector-valued continuous
functions is a generalized-lush space whenever $E$ is, and show the
stability of generalized-lush spaces by $c_0$-, $l_1$- and
$l_\infty$-sums. Then we prove that the Mazur-Ulam property holds
for a larger class of Banach spaces than GL-spaces, called
local-GL-spaces, including all lush spaces and GL-spaces.
Furthermore, we show that the $C(K,E)$ is a local-GL-space whenever
$E$ is, and the stability by $c_0$-, $l_1$- and $l_\infty$-sums also
holds for local-GL-spaces. 

Throughout this paper, we consider the spaces all over the real
field. For a Banach space $E$,  $B_E$, $S_E$, $E^*$ and $L(E)$  will stand for the unit ball of $E$, the unit sphere of $E$,
the dual space and the Banach algebra of all bounded linear operators on $E$.  A slice is a subset of $B_E$ of the form
\begin{align*}
S(x^*,\alpha)=\{x\in B_E: Re \, x^*(x)>1-\alpha\},
\end{align*}
where $x^*\in S_{E^*}$ and $0<\alpha<1$.

We recall here some basic concepts.
 \begin{definition}Let $E$ be a
Banach space.
\begin{enumerate}
\item $E$ is said to be a \emph{CL-space} if for every maximal convex set $C$ of $S_E$, we have $B_E=co(C\cup-C)$.
\item  $E$ is said to be an \emph{almost-CL-space} if for every maximal convex set $C$ of $S_E$, we have $B_E=\overline{co}(C\cup-C)$.
\item $E$ is said to be \emph{lush} if  for every $x, y \in S_E$ and every $\varepsilon > 0$,
there exists a slice $S=S(x^*, \varepsilon)$ such that $x\in S$ and dist$(y,\mbox{aco}(S))<\varepsilon$.
\end{enumerate}\end{definition}
It is an evident implication that $(1)\Rightarrow (2) \Rightarrow
(3)$, and none of the one-way implications can be reversed (see
\cite[Proposition 1]{M} and \cite[Example 3.4]{BK}).

The numerical index of a Banach space $E$ was first suggested by G.
Lumer in 1968 (see \cite{DM}), and it is the constant $n(E)$ defined
by
\begin{align*}
n(E)&=\inf\{v(T): T\in L(E),\|T\|=1\}\\
&=\max\{k\geq 0: k\|T\|\leq v(T)\  \    T\in L(E)\},
\end{align*}
where $v(T)$ is the numerical radius of $T$ and is given by
\begin{align*}
v(T)=\sup\{|x^*(T(x))|: x\in S_E, x^*\in S_{E^*},x^*(x)=1\}.
\end{align*}
More information and background on numerical indices can be found in
the recent survey \cite{BK4} and references therein.

\section{Generalized-lush spaces}

The aim of this section is to study generalized-lush spaces (GL-spaces for short). We present many examples and prove a stronger property for separable GL-spaces, and we also show that GL-spaces have some stability properties.
\begin{definition}
A Banach space $E$ is said to be a \emph{generalized-lush space \rm {(GL-space)}} if
for every $x\in S_E$ and every $\varepsilon> 0$ there exists a slice $S = S(x^*,\varepsilon)$ with $x^*\in S_{E^*}$
such that $$x\in S \quad \mbox{and }\quad \mbox{dist}(y, S)+\mbox{dist}(y,-S)<2+\varepsilon$$ for all $y\in S_E$.
\end{definition}


The following proposition for separable GL-spaces is based on an
idea from \cite[Lemma 4.2]{KM}, and it is of independent interest.

\begin{proposition} Let $E$ be a separable GL-space and $G\subset S(E^*)$ a norming subset for $E$. Then for every $\varepsilon>0$ the set
\begin{align*}
\{x^*\in G: \mbox{dist}(y, S)+\mbox{dist}(y,-S)<2+\varepsilon \
\mbox{ for all}\  y\in S_E , \mbox{ where } S=S(x^*, \varepsilon )\}
\end{align*}
is a weak$^*$ $G_\delta$-dense subset of the weak$^*$ closure of $G$.
\end{proposition}
\begin{proof}
Let $(y_n)\subset S_E $ be a sequence dense in $S_E$. Fix $0<\varepsilon<1$.  Given $n\geq1$ , set
\begin{align*}
K_n=\{x^*\in G:  \ \ \mbox{dist}(y_n, S)+\mbox{dist}(y_n,-S)<2+\varepsilon\ \ \mbox{where} \  \ S=S(x^*, \varepsilon )\}.
\end{align*}
Then $K_n$ is weak$^*$-open and $\overline{K_n}^{\omega^*}=\overline{G}^{\omega^*}$.
Indeed, if $x^*\in K_n$, then there exist $x_n\in S(x^*,\varepsilon)$ and $z_n\in -S(x^*,\varepsilon)$ such that
\begin{align*}
\|x_n-y_n\|+\|y_n-z_n\|<2+\varepsilon.
\end{align*}
Let $$U=\{y^*\in G: y^*(x_n)>1-\varepsilon \  \ \mbox {and}  \ \ y^*(-z_n)>1-\varepsilon\} .$$ Then it is easily checked that $U$ is a weak$^*$-neighborhood of $x^*$ in $G$ satisfying  $U\subset K_n$. Thus $K_n$ is weak$^*$-open.

 To prove $\overline{K_n}^{\omega^*}=\overline{G}^{\omega^*}$, it is enough to show that $G\subset \overline{K_n}^{\omega^*}$. Since \cite [Lemma 3.40]{F} states that for every $x^*\in G$, the weak$^*$-slices containing $x^*$ form a neighborhood base of $x^*$, it suffices to prove that the  weak$^*$-slice $S(x,\varepsilon_1)\cap K_n\neq\emptyset$ for all $\varepsilon_1\in (0, \varepsilon)$.  Since $E$ is a GL-space, there is a slice $S=S(y^*,\varepsilon_1/3)$ such that
\begin{align*}
x\in S\    \ \mbox{and}\   \ \mbox{dist}(y_n,S)+\mbox{dist}(y_n,-S)<2+\varepsilon_1.
\end{align*}
Thus we may find $x_n'\in S$ and $z_n' \in -S$ such that
\begin{align*}
\|x_n'-y_n\|+\|y_n-z_n'\|<2+\varepsilon_1 \  \ \mbox{and}\ \ \|x+x_n'-z_n'\|>3-\varepsilon_1.
\end{align*}
Note that $G$ is a norming subset of $S_{E^*}$. Thus there is a $z^*\in  G$ such that
\begin{align*}
z^*(x+x_n'-z_n')=\|x+x_n'-z_n'\|>3-\varepsilon_1.
\end{align*}
This implies that $z^*\in S(x, \varepsilon_1)\cap K_n$.

Now set $K=\bigcap_{n\in\mathbb{N}} K_n$. Then by the Baire theorem, $K$ is a weak$^* $ $G_\delta$-dense subset of $\overline{G}^{\omega^*}$. This together with density of $(y_n)$ in $S_E$ gives the desired conclusion.
\end{proof}

As a consequence, we have a stronger characterization for separable GL-spaces which indicates that the $x^*$ in the definition of GL-spaces can be chosen from ext$(B_{E^*})$.
\begin{corollary}\label{cor1}
Let $E$ be a separable Banach space. Then $E$ is a GL-space if and only if for every $x\in S_E$ and every $\varepsilon> 0$
there exists a slice $S = S(x^*,\varepsilon)$ with $x^*\in$ ext$(B_{E^*})$
such that $$x\in S \quad \mbox{and}\quad \mbox{dist}(y, S)+\mbox{dist}(y,-S)<2+\varepsilon$$ for all $y\in B_E$.
\end{corollary}

Now we have the following important examples.
\begin{example}\label{exa2}
Every almost-CL-space is a GL-space.
\end{example}
\begin{proof}
Let $E$ be an almost-CL-space. For every $x\in S_E$ and $\varepsilon>0$, there exists a maximal convex set $C$ of $S_E$ such that $x\in C$.
Choose $f\in S_{E^*}$ such that $f(z)=1$ for every $z\in C$, and set $S=S(f, \varepsilon)$. Then $C\subset S$.  Since $E$ is an almost-CL-space, it follows that
$B_E=\overline{co}(S\cup-S)$. So for every $y\in S_E$, there are $\lambda\in[0,1]$, $y_1\in S$ and $y_2\in -S$ such that
\begin{align*}
\|\lambda y_1+(1-\lambda)y_2-y\|<\varepsilon/2.
\end{align*}
This leads to
\begin{align*}
\| y_1-y\|+\| y_2-y\|<2+\varepsilon,
\end{align*}
which completes the proof.
\end{proof}

Since all $C(K)$, real $L_1(\mu)$ are CL-spaces (in particular, almost-CL-spaces), they are GL-spaces. Moreover, according to \cite[Theorem 4.3]{KM} showing that the separable lush space enjoys a stronger property, we can have a larger class of spaces which are GL-spaces, and they are not almost-CL-spaces in general (see, \cite [Example 3.4]{BK}).
\begin{example}
Every separable lush space is a GL-space.
\end{example}
\begin{proof}
Note that \cite[Theorem 4.3]{KM} implies that if $E$ is a separable lush space, then there is a norming subset $K$ of $S_{E^*}$ such that
\begin{align*}
B_E=\overline{co}(S(x^*, \varepsilon) \cup-S(x^*, \varepsilon))
\end{align*}
for every $x^*\in K$ and every $\varepsilon>0$. A similar analysis as in Example \ref{exa2} yields the desired conclusion.
\end{proof}

Let $K$ be a compact Hausdorff space. A closed subspace $X$ of
$C(K)$ is said to be \emph{C-rich} if for every nonempty open subset
$U$ of $K$ and every $\varepsilon>0$, there is a positive function
$h$ with $\|h\|=1$ and supp$(h)\subset U$ such that \mbox{dist}$(h,
X)<\varepsilon$. This definition covers all finite-codimensional
subspaces of $C[0,1]$ (see \cite[Proposition 2.5]{BK}) and a
subspace $X$ of $C[0,1]$ whenever $C[0,1]/X$ does not contain a copy
of $C[0,1]$ (see \cite[Proposition 1.2 and Definition 2.1]{KP}). For
more examples and results about C-rich subspaces we refer to
\cite{BK2,BK3,KM} and references therein. Notice that all C-rich
subspaces of $C(K)$ have been proved in \cite[Theorem 2.4]{BK} to be
lush. Therefore we get the following example.
\begin{example}
Every C-rich separable subspace of $C(K)$ is a GL-space.
\end{example}
Observe that all the above examples of GL-spaces are Banach spaces with numerical index 1.  We remark from the following examples that there may exist many GL-spaces whose numerical index are not 1. The two-dimensional space with hexagonal norm is firstly introduced by M. Mart\'{i}n and J. Meri \cite{MR3}.

\begin{example}\label{exa1}
The space $E=(\mathbb{R}^2,\|\cdot \|)$ whose norm is given by
\begin{eqnarray*}
\|(\xi,\eta)\|=\max\{|\eta|, |\xi|+1/2|\eta|\} \quad \forall\, (\xi,\eta)\in E,
\end{eqnarray*}
with numerical index 1/2 is a GL-space.
\end{example}
\begin{proof}
It is shown by \cite[Theorem 1]{MR3} that $E$ has numerical index 1/2. To prove that $E$ is a GL-space, given $x=(a, b)\in S_E$ and $\varepsilon>0$, we divide the proof into two cases. By symmetry considerations, we assume that $a,b\geq 0$.

Case 1: $b=1$. Define a functional $f\in S_{E^*}$ by $f(z)=\eta$ for all $z=(\xi, \eta)\in E$. Set $S=S(f, \varepsilon)$. Then $x\in S$, and for every $y=(c, d)\in S_E$, consider the two vectors
\begin{align*}
y_1=(c , 1) \ \mbox{and}  \  y_2=(c,  -1).
\end{align*}
We clearly have $y_1\in S$ and $y_2\in -S$, and moreover,
\begin{align*}
\|y-y_1\|+\|y-y_2\|=2<2+\varepsilon.
\end{align*}

Case 2: $b<1$. We make the convention that sign$(0)=1$. Let $f\in S_{E^*}$ be defined by $f(z)=\xi+ \eta/2$ for every $z=(\xi, \eta)\in E$.
This guarantees that $x\in S=S(f, \varepsilon)$. For every $y=(c, d)\in S_E$, we set
\begin{equation*}
\begin{cases}
y_1=(\mbox{sign}(c), 0), \  y_2=\mbox{sign}(d)(1/2, 1) \hspace{1.3cm}\text{if } \hspace{0.1cm} cd\leq0; \\
y_1=-(\mbox{sign}(c), 0),\  y_2=\mbox{sign}(d)(1/2, 1) \hspace{1cm}\text{if } \hspace{0.1cm} cd>0 \ \mbox{and} \ |d|=1; \\
y_1=y,\ y_2=-y,                                        \hspace{4.9cm}\text{if } \hspace{0.1cm}cd>0 \ \mbox{and} \ |d|<1.

\end{cases}
\end{equation*}
 Then $y_1, y_2\in S\cup(-S)$ satisfy
\begin{align*}
\|y-y_1\|+\|y-y_2\|=2<2+\varepsilon.
\end{align*}
We thus complete the proof.
\end{proof}

By Example \ref{exa1}, the following Theorems \ref{lem5},
\ref{prop1} and \cite[Proposition 1]{M2} which shows that the numerical index of the $c_0$-,$l_1$-,or $l_\infty$-sum of Banach spaces is the infimum numerical index of the summands, we may construct more
examples of specific GL-spaces with numerical index 1/2.
\begin{example}
The space $E=(c_0,\|\cdot\|)$ equipped with the norm
\begin{align*}
\|x\|=\max\{\sup_{k\in \mathbb{N}}|\xi_k|, |\xi_1|+1/2|\xi_2|\}
\quad \forall x=(\xi_k)\in E
\end{align*}
is a GL-space with numerical index 1/2.
\end{example}
\begin{proof} It is actually the space $c_0\bigoplus_\infty X$ where $X$ is
just the hexagonal space in Example \ref{exa1}.
\end{proof}

We shall give an observation that in the definition of GL-spaces we
can take $y$ to be in the unit ball instead of being in the unit
sphere. With the help of this observation, one can check whether the
space being considered is a GL-space in an easier way. We will use it
later to get some stability properties of GL-spaces.
\begin{lemma}\label{rem2}
 If $E$ is a GL-space, then for every $x\in S_E$ and every $\varepsilon> 0$ there exists a slice $S = S(x^*,\varepsilon)$ with $x^*\in S_{E^*}$
such that $$x\in S \quad \mbox{and }\quad \mbox{dist}(y,
S)+\mbox{dist}(y,-S)<2+\varepsilon$$ for all $y\in B_E$.
\end{lemma}
\begin{proof}For every $x\in S_E$ and every $\varepsilon> 0$, let $S = S(x^*,\varepsilon)$ be such that
$$x\in S \quad \mbox{and }\quad \mbox{dist}(z, S)+\mbox{dist}(z,-S)<2+\varepsilon$$ for all $z\in S_E$.
Given $y\in B_E$, since the case where $y=0$ is trivial, we may
assume that $y\neq0$. Then there exist $u, -v\in S$ such that
$$\|u-\frac{y}{\|y\|}\|+\|v-\frac{y}{\|y\|}\|<2+\varepsilon.$$
Triangle inequality hence yields
\begin{align*}
\|u-y\|+\|v-y\|<2+\varepsilon\|y\|\leq 2+\varepsilon.
\end{align*}
This completes the proof.
\end{proof}

Given a compact Hausdorff space $K$ and a Banach space
$E$, we denote by $C(K,E)$ the Banach space of all continuous functions from $K$ into
$E$, endowed with its natural supremum norm.
\begin{theorem}\label{lem5}
Let $K$ be a compact Hausdorff space and $E$ a GL-space, then $C(K,E)$ is a GL-space.
\end{theorem}
\begin{proof}
Given $f\in S_{C(K,E)}$ and $\varepsilon>0$, there exists a $t_0\in K$ such that $\|f(t_0)\|=1$. Since $E$ is a GL-space, it follows from this and Lemma \ref{rem2} that there is an $x^*\in S_E^*$ with $S_{x^*}=S(x^*, \varepsilon/2)$ such that
\begin{align*}
f(t_0)\in S_{x^*}\  \  \mbox{and} \  \  \mbox{dist}(y, S_{x^*})+\mbox{dist}(y,-S_{x^*})<2+\varepsilon/2
\end{align*}
for all $y\in B_E$.
Define a functional $f^*\in S_{{C(K,E)}^*}$ by $f^*(g)=x^*(g(t_0))$ for every $g\in {C(K,E)}$, and put $S=S(f^*, \varepsilon)$.
 For every $g\in S_{C(K,E)}$, we have $g(t_0)\in B_E$. Thus there are $y_1 \in S_{x^*}$ and $y_2\in-S_{x^*}$ such that
 \begin{align*}
 \|g(t_0)- y_1\|+\|g(t_0)-y_2\|<2+\varepsilon/2.
 \end{align*}
Then we can build a continuous map $\phi: K\rightarrow [0,1]$ defined by
\begin{align*}
\phi(t_0)=1\  \ \mbox{and} \ \  \phi(t)=0  \ \    \mbox{if}  \   \   \|g(t)-g(t_0)\|\geq\varepsilon/4 .
\end{align*}
Consider $h_1\in S $ and  $h_2\in -S$ given by
\begin{align*}
h_i(t)=\phi(t)y_i+(1-\phi(t))g(t)(i=1, 2) \ \  \mbox{for every}\ t\in K.
\end{align*}
Then it is trivial to see that
\begin{align*}
\|g-h_1\|+\|g-h_2\|<2+\varepsilon.
\end{align*}
Hence $C(K,E)$ is a GL-space.
\end{proof}

For more examples of GL-spaces, we need discuss the stability of GL-spaces
by $c_0$-, $l_1$- and $l_\infty$-sums. Recall that the  $c_0$-sum (resp.  $l_1$-sum and $l_\infty$-sum) of a family of Banach spaces $\{E_\lambda : \lambda\in\Lambda\} $ are denoted by $[\bigoplus_{\lambda\in\Lambda} E_\lambda]_{c_0}$ (resp. $[\bigoplus_{\lambda\in\Lambda} E_\lambda ]_{l_1}$ and $[\bigoplus_{\lambda\in\Lambda} E_\lambda ]_{l_\infty}$). 
\begin{theorem}\label{prop1}
Let $\{E_\lambda : \lambda\in\Lambda\}$ be a family of Banach spaces, and let
$E = [\bigoplus_{\lambda\in\Lambda} E_\lambda]_{F}$ where $F=c_0 \,  \mbox{or } \,  l_1$.
Then $E$  is a GL-space if and only if each $E_\lambda$ is a GL-space.
\end{theorem}
\begin{proof}
Note that $E^*=[\bigoplus_{\lambda\in\Lambda} E_\lambda^*]_{l_1}$ if $F=c_0$ and  $E^*=[\bigoplus_{\lambda\in\Lambda} E_\lambda^*]_{l_\infty}$ if $F=l_1$. This fact will
be used without comment in the following proof.

In the $c_0$-sum case, we first show the ``if" part. Fix $x=(x_\lambda)\in S_E$ and $\varepsilon>0$. We may find a $\lambda_0$ such that $\|x_{\lambda_0}\|=1$.  Since $E_{\lambda_0}$ is a GL-space, by Lemma \ref{rem2} there is a slice $S_{\lambda_0}=S(x_{\lambda_0}^*,\varepsilon) \subset B_{E_{\lambda_0}}$ with $x_{\lambda_0}^*\in S_{E_{\lambda_0}^*}$ such that $$ x_{\lambda_0}\in S_{\lambda_0}\ \  \mbox{and} \   \ \mbox{dist}(z, S_{\lambda_0})+\mbox{dist}(z,-S_{\lambda_0})<2+\varepsilon$$ for all $z\in B_{E_{\lambda_0}}$. Choose $x^*=(x_\lambda^*)\in S_{E^*}$ with $x_\lambda^*=0$ for all $\lambda\neq \lambda_0$, and let $S=S(x^*, \varepsilon).$
 Then $x\in S$, and it is easy to see from the definition of $E$ that
\begin{align}\label{equ:110}
\mbox{dist}(y, S)+\mbox{dist}(y,-S)<2+\varepsilon
\end{align}
for all $y\in S_E$. Thus $E$ is a GL-space.

 Now we deal with the ``only if" part. For every $\lambda\in\Lambda$, fix $x_\lambda \in S_{E_\lambda}$ and $ \varepsilon>0$. Take $x=(x_\delta)\in S_E$ with  $x_\delta=0$ for all $\delta\neq \lambda$. Then $x\in S_E,$ and thus there exists an $x^*=(x_\delta^*)\in S_{E^*}$ with $S=S(x^*, \varepsilon/2)$
 such that
 \begin{align}\label{equ:20}
  x\in S\ \  \mbox{and} \   \ \mbox{dist}(y, S)+\mbox{dist}(y,-S)<2+\varepsilon/2
\end{align}
for all $y\in S_E$. Note that $x_\lambda\in S_\lambda=S(x_\lambda^*/\|x_\lambda^*\|,\varepsilon)$. To show that $E_\lambda$  is a GL-space, it remains to check that for all $y_\lambda\in S_{E_\lambda}$
 \begin{align*}
 \mbox{dist}(y_\lambda, S_\lambda)+\mbox{dist}(y_\lambda,-S_\lambda)<2+\varepsilon.
\end{align*}

Now given $y_\lambda\in S_{E_\lambda}$, consider $y=(y_\delta)\in S_E$ with $y_\delta =0$ for all $\delta\neq \lambda$. By (\ref{equ:20}), there are  $u=(u_\delta)\in S$ and $v=(v_\delta)\in -S$ such that
\begin{align*}
\|y-u\|+\|y-v\|<2+\varepsilon/2.
\end{align*}
The definition of $E$ thus gives
\begin{align*}
\|y_\lambda-u_\lambda\|+\|y_\lambda-v_\lambda\|<2+\varepsilon/2.
\end{align*}
Observe that $\|x_\lambda^*\|\geq x_\lambda^*(x_\lambda)>1-\varepsilon/2$, and therefore $\sum_{\delta\neq \lambda}\|x_\delta^*\|<\varepsilon/2$.
So
\begin{align*}
x_\lambda^*(u_\lambda)>1-\varepsilon/2-\sum_{\delta\neq \lambda}\|x_\delta^*\|>1-\varepsilon .
\end{align*}
Similarly, $ x_\lambda^*(-v_\lambda)>1-\varepsilon.$
   Hence $E_\lambda$  is a GL-space.

In the $l_1$-sum case, let us prove the ``if" part. Given $x=(x_\lambda)\in S_E$ and $\varepsilon>0$, for each $\lambda$ with $x_\lambda\neq 0$, there is a corresponding slice  $S_\lambda=S(x_\lambda^*, \varepsilon)$ with $x_\lambda^*\in S_{E_\lambda^*}$ such that
\begin{align*}
x_\lambda^*(x_\lambda)>(1-\varepsilon )\|x_\lambda\| \  \mbox{and} \   \ \mbox{dist}(z_\lambda, S_\lambda)+\mbox{dist}(z_\lambda,-S_\lambda)<2+\varepsilon
\end{align*}
for all $z_\lambda\in S_{E_\lambda}$. Then $x^*=(x_\lambda^*)\in S_{E^*}$ with $x_\lambda^*=0$ whenever $x_\lambda=0$, and the required slice satisfying (\ref{equ:110}) is $S(x^*,\varepsilon)$. Therefore $E$  is a GL-space.

For the ``only if" part, fix $x_\lambda\in S_{E_\lambda}$ and $0<\varepsilon<1/2$. Then $x=(x_\delta)\in S_E$ where $x_\delta =0$ for all $\delta\neq \lambda$. Since $E$  is a GL-space, there is an $x^*=(x_\delta^*)\in S_{E^*}$ with $S=S(x^*,\varepsilon/4)$ such that
$$ x\in S\ \  \mbox{and} \   \ \mbox{dist}(y, S)+\mbox{dist}(y,-S)<2+\varepsilon/4$$ for all $y\in S_E$.  We shall prove that the slice $S_\lambda=S(x_\lambda^*/\|x_\lambda^*\| , \varepsilon)$ is the desired one, namely that $x_\lambda\in S_\lambda$ and   $\mbox{dist}(y_\lambda, S_\lambda)+\mbox{dist}(y_\lambda,-S_\lambda)<2+\varepsilon$ for all $y_\lambda\in S_{E_\lambda}$.

It is easily checked that $x_\lambda\in S_\lambda$. For every $y_\lambda\in S_{E_\lambda}$, since $y=(y_\delta)$ is in $S_E$ where $y_\delta =0$ for all $\delta\neq \lambda$,  there are $u=(u_\delta)\in S$ and $v=(v_\delta)\in -S$ such that
\begin{align}\label{equ:10}
\|y-u\|+\|y-v\|<2+\varepsilon/4.
\end{align}
It follows from the definition of $E$ that
\begin{align}
\|y-u\|+\|y-v\|&=\|y_\lambda-u_\lambda\|+\sum_{\delta\neq \lambda}\|u_\delta\|+\|y_\lambda-v_\lambda\|+\sum_{\delta\neq \lambda}\|v_\delta\| \nonumber\\
&>\|y_\lambda-u_\lambda\|+1-\varepsilon/4-\|u_\lambda\|+\|y_\lambda-v_\lambda\|+1-\varepsilon/4-\|v_\lambda\|\nonumber \\
&=\|y_\lambda-u_\lambda\|-\|u_\lambda\|+\|y_\lambda-v_\lambda\|-\|v_\lambda\|+2-\varepsilon/2.\label{equ:11}
\end{align}
We deduce from (\ref{equ:10}) and (\ref{equ:11}) that
\begin{align*}
\|u_\lambda\|>1/2-\varepsilon/2 \ \ \mbox{and}\ \  \|v_\lambda\|>1/2-\varepsilon/2.
\end{align*}
Hence
\begin{align*}
x_\lambda^*(u_\lambda)>1-\varepsilon/4-\sum_{\delta\neq \lambda}\|u_\delta\|\geq1-\varepsilon/4-1+\|u_\lambda\|\geq(1-\varepsilon)\|u_\lambda\|,
\end{align*}
and similarly,
\begin{align*}
x_\lambda^*(-v_\lambda)>(1-\varepsilon)\|v_\lambda\|.
\end{align*}

So $w_\lambda=u_\lambda/\|u_\lambda\|$ and $t_\lambda=-v_\lambda/\|v_\lambda\|$ are in $S_\lambda$. The desired estimate
\begin{align*}
\|y_\lambda-w_\lambda\|+\|y_\lambda+t_\lambda\|<2+\varepsilon
\end{align*}
which follows from (\ref{equ:11}) completes the proof.
\end{proof}

We also have a proposition establishing that the class of GL-spaces
is stable under the  $l_\infty$-sum, and we omit the proof since it
is just a slight modification of the ``if" part in the $c_0$-case.
\begin{proposition}\label{infty}
Let $\{E_\lambda : \lambda\in\Lambda\}$ be a family of GL-spaces, and let
$E = [\bigoplus_{\lambda\in\Lambda} E_\lambda]_{l_\infty}$.
Then $E$ is a GL-space.
\end{proposition}

\section{The Mazur-Ulam property for local-GL-spaces} \label{se2}

The main aim of this section is to prove that a larger class of
Banach spaces have the Mazur-Ulam property. We begin with a proposition which is the key step to prove Theorem \ref{mth}.
\begin{proposition} \label{lem1}
Let $E, F$ be Banach spaces, and let $T:S_E\to S_F$ be an isometry
(not necessarily surjective). If $E$ is a GL-space, then we have
\begin{equation*}
    \|T(x)-\lambda\, T(y)\|\ge\|x-\lambda\, y\| \ \ \mbox{ for all } x,y\in S_E \mbox{ and } \lambda\ge 0.
\end{equation*}
\end{proposition}
\begin{proof}
Given $x,y\in S_E$ with $x\neq y$ and $\lambda>0$, set
\begin{align*}
z=\frac{x-\lambda y}{\|x-\lambda y\|}.
\end{align*}
Since $E$ is a GL-space, given $\varepsilon> 0$, there exists a
functional $f\in S_{E^*}$ with $S=S(f, \varepsilon)$ such that
\begin{align*}
z\in S \quad \mbox{and} \quad  \mbox{dist}(w,S)+\mbox{dist}(w,-S)<2+\varepsilon
\end{align*}
for all $w\in S_E$. Therefore, there exist $x_1, y_1\in S$  and $x_2,y_2\in -S$ such that
\begin{equation*}
    \|x-x_1\|+\|x-x_2\|<2+\varepsilon \ \ \mbox{ and } \ \ \|y-y_1\|+\|y-y_2\|<2+\varepsilon.
\end{equation*}
Then
\begin{align*}
2-2\varepsilon <f(x_1)-f(x)+f(x)-f(x_2) \leq
\|x-x_1\|+\|x-x_2\|<2+\varepsilon.
\end{align*}
This implies that
\begin{align}\label{equ:1}
f(x_1)-f(x) \ge  \|x-x_1\| - 3\varepsilon.
\end{align}
A similar analysis gives
\begin{align}\label{equ:2}
f(y)-f(y_2) \ge  \|y-y_2\| - 3\varepsilon.
\end{align}

For $i=1$ or 2, replace $x_i $ by $x_i/\|x_i\|$  and  $y_i$ by $y_i/\|y_i\|$ respectively if necessary we may assume that $x_i$ and $y_i$ have norm 1.  Then there exists a functional $g \in S_{F^*}$ such that
\begin{align*}
g(T(x_1))-g(T(y_2))=\|T(x_1)-T(y_2)\|=\|x_1-y_2\|> 2-2\varepsilon.
\end{align*}
It follows that
\begin{align*}
g(T(x_1)) > 1-2\varepsilon \quad \mbox{and} \quad g(T(y_2)) < -1+2\varepsilon.
\end{align*}
Thus by (\ref{equ:1}) and (\ref{equ:2}), we have
 \begin{align*}
f(x)&\le f(x_1)-\|x-x_1\|+3\varepsilon
 \\  &\le 1-\|T(x)-T(x_1)\|+3\varepsilon
 \\   &\le 1-(g(T(x_1))-g(T(x))+3\varepsilon
   \\   &\le g(T(x))+5\varepsilon
\end{align*}
and
\begin{align*}
f(y)&\ge f(y_2)+\|y-y_2\|-3\varepsilon
 \\  &\ge -1+\|T(y)-T(y_2)\|-3\varepsilon
 \\   &\ge-1+(g(T(y))-g(T(y_2))-3\varepsilon
   \\   &\ge g(T(y))-5\varepsilon.
\end{align*}
As a consequence,
\begin{align*}
\|x-\lambda y\|(1-\varepsilon)< f(x-\lambda y)&\leq g(T(x))+5\varepsilon-\lambda g(T(y))+5\lambda \varepsilon  \\
                                                 & \le\|T(x)-\lambda T(y)\|+ (5+5\lambda)\varepsilon.
\end{align*}
Since $\varepsilon$ can be arbitrarily small, we complete the proof.
\end{proof}

\begin{theorem}\label{th1}
Every GL-space $E$ has the MUP.
\end{theorem}
\begin{proof}
Let $F$ be a Banach space, and let $T: S_E\rightarrow S_F$ be a surjective isometry. We need to show that $T$
can be extended to a linear surjective isometry from $E$ onto $F$.
We first claim that for all $x,y\in S_E$ and $\lambda\ge 0$.
\begin{align}\label{equ:4}
\|T(x)-\lambda T(y)\|= \|x-\lambda y\|.
\end{align}
Otherwise by Proposition \ref{lem1}, there exist $\lambda_0>0, x_0, y_0 \in S_E$ such that
\begin{align}\label{equ:3}
\|T(x_0)-\lambda_0 T(y_0)\|>\|x_0-\lambda_0 y_0\|.
\end{align}
Replace $\lambda_0$ by $1/\lambda_0$ if necessary we may assume that $\lambda_0<1$. Since $\|\lambda_0 T(y_0)\|=\lambda_0<1$,
there exists $T(v) \in S_F$ with $v\in S_E$ such that $\lambda_0 T(y_0)$ belongs to the segment $(T(x_0), T(v))$ of $B_F$.
By (\ref{equ:3}) and Proposition \ref{lem1} we have
\begin{align*}
\|v-x_0\|=\|T(v)-T(x_0)\|&=\|T(v)-\lambda_0 T(y_0)\|+\|\lambda_0
T(y_0)-T(x_0)\|\\ &>\|v-\lambda_0 y_0\|+\|\lambda_0 y_0-x_0\|\\
&\geq \|v-x_0\|.
\end{align*}
It is a contradiction.
Now we may define the required extension $\widetilde{T}$ of $T$ by
$$\widetilde{T}(x)= \left
\{ \begin{array}{ll}
\|x\|T(\frac{x}{\|x\|}), & \mbox{ if }\;x\neq0;\\
0, & \mbox{ if }\;x=0.
\end{array}
\right.$$
It is easily seen from (\ref{equ:4}) that $\widetilde{T}: E \rightarrow F$ is a surjective isometry whose restriction to the unit sphere $S_E$ is just
$T$. The Mazur-Ulam theorem hence shows that $\widetilde{T}$ is linear as desired. The proof is complete.
\end{proof}
Note that the technique in the proof of Theorem \ref{th1} is still valid in more general case. We now state a result here since it will be of use later.
\begin{proposition}\label{th2}
Let $E, F$ be Banach spaces, and let
$T:S_E\to S_F$ be a surjective isometry such that
\begin{equation*}
    \|T(x)-\lambda\, T(y)\|\ge\|x-\lambda\, y\| \ \ \mbox{ for all } x,y\in S_E \mbox{ and } \lambda\ge 0.
\end{equation*}
Then $E$ has the MUP.
\end{proposition}

Now we introduce a class of spaces called local-GL-spaces (including
GL-spaces and lush spaces) which have the MUP. This definition is a
weakening of the notion of lush spaces in the real case. We can see
from the above Example \ref{exa1} that this weakening is strict.

\begin{definition} A Banach space $E$ is said to be a \emph{local-GL-space} if for every separable subspace $X\subset E $,
there is a GL-subspace $Y\subset E$ such that $X \subset Y \subset
E$.
\end{definition}
\begin{example} GL-spaces are local-GL-spaces.
\end{example}
The equivalent definition of lush space \cite[Theorem
4.2]{BK2} proves the following.
\begin{example}
Lush spaces are local-GL-spaces.
\end{example}


We now present the main result of this section.
\begin{theorem}\label{mth}
Every local-GL-space has the MUP.
\end{theorem}
\begin{proof}
Let $E$ be a local-GL-space, $F$ a Banach space and $T:
S_E\rightarrow S_F$ a surjective isometry. We next show that $T$
can be extended to a linear surjective isometry from $E$ onto $F$.

Fix $x, y\in S_E$. Let $X=\mbox{span}(x,y)$.
Since $E$ is a local-GL-space, there is a GL-space $Y\subset E$
such that $X\subset Y$. We consider $T$ to be an isometry from $S_Y$ to
$S_F$. Then Propositions \ref{lem1} and \ref{th2} clearly lead to
the fact that $T$ can be extended to a linear surjective isometry
from $E$ onto $F$.
\end{proof}
We emphasize two evident consequences of the above theorem.
\begin{corollary}
Every lush space has the MUP.
\end{corollary}

\begin{corollary}
Every C-rich subspace of $C(K)$ has the MUP.
\end{corollary}

By the following properties, we can get more examples of spaces
having the MUP.

\begin{proposition}\label{prop5}
If $E$ is a local-GL-space, then $C(K,E)$ is a local-GL-space.
\end{proposition}
\begin{proof}
Let $X$ be a separable subspace of $C(K,E)$. We shall prove that the
set $$E_X=\bigcup_{t\in K} \{f(t): f\in X\}$$ is a separable subset
of $E$. Let $\{f_n\}$ be a dense sequence of $X$. Given $n, m\geq 1$
and $s\in K$, set
$
V_{s,\, m,\, n}=\{t\in K:  \|f_n(t)-f_n(s)\|<1/m\}.
$
The compactness of $K$ implies that there is a finite subset
$\{s_{i}^{m,n}: i=1,\cdots, k_{m,n} \}$ of $K$ such that
$K=\bigcup_{i=1}^{k_{m,n}} V_{s_{i}^{m,n},\, m,\,n}$. Then it is an
elementary check that the set
\begin{align*}
M=\bigcup _{n=1}^\infty \bigcup_{m=1}^\infty \{f_n(s_i^{m,n}):
i=1,\cdots k_{m,n} \}
\end{align*}
is a dense subset of $E_X$. It follows that
$N_X=\overline{\mbox{span} \{E_X\}}$ is a separable subspace of $E$.
Note that the $E$ is a local-GL-space. So we may find a
GL-space $M_X$ such that $N_X\subset M_X \subset
E$.

Let $Y=C(K, M_X)$. Then $X\subset Y $, and Theorem \ref{lem5} shows
that $Y$ is a GL-space.  This completes the proof.
\end{proof}

\begin{corollary}
Let $E$ be a local-GL-space and $K$ be a compact Hausdorff space.
Then $C(K,E)$ has the MUP.
\end{corollary}

The proof of Theorem \ref{prop1} can be adapted to yield a characterization of the $c_0$-, $l_1$-sums of lush spaces in both real and complex cases.
The ``if" part of it has been noted in \cite[Propsosition 5.3]{BK2}, and the ``only if" part is probably known but we include an argument here (as
we do not find it explicitly stated in the literature).
\begin{proposition}\label{prop2}
Let $\{E_\lambda : \lambda\in\Lambda\}$ be a family of Banach spaces, and let
$E = [\bigoplus_{\lambda\in\Lambda} E_\lambda]_{F}$ where $F=c_0 \,  \mbox{or} \,  l_1$. Then $E$ is a lush space if and only if $E_\lambda$ is a  lush space for every $\lambda\in\Lambda$.
\end{proposition}
\begin{proof}
It has been proved in \cite[Propsosition 5.3]{BK2} that each $E_\lambda$ is a lush space, then $E$ is also lush. We only check the ``only if" statement. Note that the $c_0$-case follows from the proof of Theorem \ref{prop1} with minor modifications. We omit the proof, leaving routine details to the readers.

Now for the $l_1$-case, fix $x_\lambda, y_\lambda \in S_{E_\lambda}$ and $0<\varepsilon<1/2$. Consider $x=(x_\delta), y=(y_\delta)\in S_E$ with $x_\delta =y_\delta=0$ for all $\delta\neq \lambda$. Then there is an $x^*=(x_\delta^*)\in S_{E^*}$ with $S=S(x^*,\varepsilon/8)$ such that
$$ x\in S\ \  \mbox{and} \   \ \mbox{dist}(y, \mbox{aco}(S))<\varepsilon/8.$$
This implies that $x_\lambda\in S_{x^*_\lambda}=S(x^*_\lambda/\|x^*_\lambda\|,\varepsilon) $ and produces a finite number of elements $\{u^i\}_{i=1}^{n}\subset S$ with $u^i=(u_\delta^i)$ and a finite number of scalars $\{\lambda_i\}_{i=1}^{n}$  with $\sum_{i=1}^{n}|\lambda_i|=1$ such that
\begin{align}\label{Ht}
\|y_\lambda-\sum_{i=1}^{n}\lambda_i u_\lambda^i\|+\sum_{\delta\neq\lambda}\|\sum_{i=1}^{n}\lambda_i u_\delta^i\|<\varepsilon/8.
\end{align}
Set $$I=\{i\in\{1,\cdots,n\}: \|u_\lambda^i\|>1/2-\varepsilon/2\}.$$ We clearly have from (\ref {Ht}) that $\|\sum_{i=1}^{n}\lambda_i u_\lambda^i\|>1-\varepsilon/8$. We then deduce from this that
$\sum_{i\in I}|\lambda_i|\geq 1-\varepsilon /4.$ The same technique in Theorem \ref{prop1} thus proves that $\widetilde{u_\lambda^i}=u_\lambda^i/\|u_\lambda^i\|\in S_{x^*_\lambda}$ for all $i\in I$, and
\begin{align}\label{Hs}
\|y_\lambda-\sum_{i\in I}\widetilde{\lambda_i }\widetilde{u_\lambda^i}\|<\varepsilon
\end{align}
where $\widetilde{\lambda_i}=\lambda_i/(\sum_{i\in I}|\lambda_i|)$. For (\ref{Hs}) we need the inequality
$$\|\sum_{i\in I}\lambda_i u_\lambda^i-\sum_{i\in I}{\lambda_i }\widetilde{u_\lambda^i}\|
\leq\sum_{i\in I}{|\lambda_i |}(1-\|u_\lambda^i\|)\leq 1-\|\sum_{i\in I}\lambda_i u_\lambda^i\|\leq 3\varepsilon/8.$$
This finishes the proof.
\end{proof}

We next give an analogue of Proposition \ref{prop2} for
local-GL-spaces. The proof of this result is routine based on
Theorem \ref{prop1}.
\begin{proposition}\label{prop3}
Let $\{E_\lambda : \lambda\in\Lambda\}$ be a family of Banach
spaces, and let $E = [\bigoplus_{\lambda\in\Lambda} E_\lambda]_{F}$
where $F=c_0 \,  \mbox{or} \,  l_1$. Then $E$ is a local-GL-space
if and only if $E_\lambda$ is a local-GL-space for every
$\lambda\in\Lambda$.
\end{proposition}

\begin{proof}
Let $P_\lambda$ be the projection of $E$ onto $E_\lambda$, and let
$I_\lambda$ be the injection of $E_\lambda$ into $E$.

 We first show
the ``if" part. Fix a separable subspace $X$ of $E$. Then
$P_\lambda(X)\subset E_\lambda$ is separable. Since $E_\lambda$ is a
local-GL-space, there is a GL-space $Y_\lambda\subset E_\lambda$
such that $P_\lambda(X)\subset Y_\lambda$. Then
$Y=[\bigoplus_{\lambda\in\Lambda} Y_\lambda]_{F}$ containing $X$ is
a subspace of $E$. Moreover it follows from Theorem \ref{prop1} that
$Y$ is a GL-space, and hence $E$ is a local-GL-space.

Now let us deal with the ``only if" part. Given $\lambda\in \Lambda$, let $X_\lambda$ be a separable subspace of $E_\lambda.$ Since $E$ is a local-GL-space,
 there is a GL-space $Y$ such that $I_\lambda(X_\lambda)\subset Y\subset E$.  Note from Theorem \ref{prop1} that $Y_\lambda=P_\lambda(Y)$ is a GL-space such that $X_\lambda\subset Y_\lambda\subset E_\lambda$.  Thus $E_\lambda$ is a local-GL-space.
\end{proof}
A similar analysis as the above proposition yields the following result.
\begin{proposition}\label{prop7}
Let $\{E_\lambda : \lambda\in\Lambda\}$ be a family of
local-GL-spaces and let $E = [\bigoplus_{\lambda\in\Lambda}
E_\lambda]_{l_\infty}$. Then $E$ is a local-GL-space.
\end{proposition}
As immediate consequences of the propositions above, we obtain that:
\begin{corollary}
Let $\{E_\lambda : \lambda\in\Lambda\}$ be a family of
local-GL-spaces. Then the space $E = [\bigoplus_\lambda
E_\lambda]_{F}, \, \mbox{where}\, F=c_0, l_1 \, \mbox{or} \ \,
l_\infty $ has the MUP.
\end{corollary}

 Throughout this paper, we can see that the
geometry properties, isometric extension, and even the numerical index
on unit spheres have harmonious inner relationship and may provide a
possible way to solve the isometric extension problem in more
general case. Note that there exist examples of Banach spaces with numerical index 1 but not lush spaces (see
\cite[Remark 4.2]{BK3}).
 Then the first natural question to ask is the
following:
\begin{question} Does every Banach space with numerical index 1 have
the MUP?
\end{question}
\noindent\textbf{Acknowledgment.} The authors express their appreciation to Guanggui Ding
for many very helpful comments regarding isometric theory in Banach spaces.

\end{document}